\documentclass{elsart3-1}
\usepackage{amsmath}
\usepackage{amssymb,latexsym}
\usepackage{amscd}
\usepackage{xspace}
\usepackage{verbatim}

\usepackage[english,francais]{babel}

\newtheorem{e-proposition}[theorem]{Proposition}

\newtheorem{e-definition}[theorem]{Definition\rm}

\renewcommand{\theequation}{\arabic{equation}}
\def\og{\leavevmode\raise.3ex\hbox{$\scriptscriptstyle\langle\!\langle$~}}
\def\fg{\leavevmode\raise.3ex\hbox{~$\!\scriptscriptstyle\,\rangle\!\rangle$}}
\newcommand{\txt}[1]{\;\text{ #1 }\;}

\newcommand{\be}{\begin{equation}}
\newcommand{\bel}[1]{\begin{equation}\label{#1}}
\newcommand{\ee}{\end{equation}}
\newtheorem{subn}{\name}

\newcommand{\bsn}[1]{\def\name{#1$\!\!$}\begin{subn}}
\newcommand{\esn}{\end{subn}}
\newtheorem{sub}{\name}[section]
\newcommand{\bs}{\begin{sub}}
\newcommand{\es}{\end{sub}}
\newcommand{\bsl}[1]{\begin{sub}\label{#1}}
\newcommand{\bth}[1]{\def\name{Theorem}\begin{sub}\label{t:#1}}
\newcommand{\blemma}[1]{\def\name{Lemma}\begin{sub}\label{l:#1}}
\newcommand{\bcor}[1]{\def\name{Corollary}\begin{sub}\label{c:#1}}
\newcommand{\bdef}[1]{\def\name{Definition}\begin{sub}\label{d:#1}}
\newcommand{\bprop}[1]{\def\name{Proposition}\begin{sub}\label{p:#1}}
\newcommand{\bnote}[1]{\def\name{\mdseries\scshape Notation}\begin{sub}\label{n:#1}}
\newcommand{\bcom}{}
\newcommand{\req}{\eqref}
\newcommand{\rth}[1]{Theorem~\ref{t:#1}}

\newcommand{\BA}{\begin{array}}
\newcommand{\EA}{\end{array}}
\newcommand{\BAN}{\renewcommand{\arraystretch}{1.2}
\setlength{\arraycolsep}{2pt}\begin{array}}
\newcommand{\BAV}[2]{\renewcommand{\arraystretch}{#1}
\setlength{\arraycolsep}{#2}\begin{array}}
\newcommand{\BSA}{\begin{subarray}}
\newcommand{\ESA}{\end{subarray}}
\newcommand{\BAL}{\begin{aligned}}
\newcommand{\EAL}{\end{aligned}}
\newcommand{\BALG}{\begin{alignat}}
\newcommand{\EALG}{\end{alignat}}
\newcommand{\BALGN}{\begin{alignat*}}
\newcommand{\EALGN}{\end{alignat*}}

\newcommand{\forevery}{\quad \forall}

\newcommand{\2}{\\[2mm]}

\newcommand{\set}[1]{\{#1\}}
\def\({{\rm (}}
\def\){{\rm )}}

\def\uar{\uparrow}
\newcommand{\abs}[1]{\left |#1\right |}
\newcommand{\norm}[1]{\left \|#1\right \|}

\newcommand{\rec}[1]{\frac{1}{#1}}
\newcommand{\opname}[1]{\mathrm{#1}\,}
\newcommand{\supp}{\opname{supp}}
\newcommand{\dist}{\opname{dist}}

\newcommand{\q}{\quad}

\newcommand{\prt}{\partial}
\newcommand{\sms}{\setminus}
\newcommand{\ems}{\emptyset}

\newcommand{\tl}{\tilde}

\newcommand{\sbs}{\subset}

\newcommand{\Consy}{Consequently\xspace}

\newcommand{\seq}{sequence\xspace}

\newcommand{\bdw}{\partial\Gw}

\newcommand{\gsmod}{$\gs$-moderate\xspace}

\def\bcom{}
\def\ga{\alpha}            
             \def\ge{\epsilon}

            \def\gl{\lambda}
\def\gm{\mu}                 
    \def\gr{\rho}        
\def\gs{\sigma}       
      
                \def\gz{\zeta}
     \def\Gd{\Delta}

\def\Gw{\Omega}              

\def\BBP {\mathbb P}   \def\BBR {\mathbb R}

\def\BBZ {\mathbb Z}

\begin{document}
\centerline{}
\begin{frontmatter}


\selectlanguage{english}
\title{Maximal solutions of equation $\Delta u=u^q$ in arbitrary domains}


\selectlanguage{english}
\author[authorlabel1]{Moshe Marcus},
\ead{marcusm@math.technion.ac.il}
\author[authorlabel2]{Laurent V\'eron}
\ead{veronl@univ-tours.fr}

\address[authorlabel1]{Department of Mathematics, Technion\\
 Haifa 32000, ISRAEL}
\address[authorlabel2]{Laboratoire de Math\'ematiques et Physique Th\'eorique, Facult\'e des Sciences\\
Parc de Grandmont, 37200 Tours, FRANCE}


\medskip

\begin{abstract}
\selectlanguage{english}
We prove bilateral capacitary estimates for the maximal solution
$U_F$ of $-\Delta u+u^q=0$ in the complement of an arbitrary
closed set $F\subset\mathbb R^N$,
 involving the Bessel capacity $C_{2,q'}$, for $q$ in the
 supercritical range $q\geq q_{c}:=N/(N-2)$.
 We derive a  pointwise necessary and sufficient
condition, via a Wiener type criterion, in order that
$U_F(x)\to\infty$ as $x\to y$ for given $y\in\prt F$.
 Finally we prove a general uniqueness
result for large solutions.  {\it To cite this article: M. Marcus,
L. V\'eron, C. R. Acad. Sci. Paris, Ser. I XXX (2007).}

\vskip 0.5\baselineskip

\selectlanguage{francais}
\noindent{\bf R\'esum\'e} \vskip 0.5\baselineskip \noindent
{\bf Solutions maximales de $\Delta u=u^q$ dans un domaine arbitraire.}
Nous d\'emontrons une estimation capacitaire bilat\'erale de la solution maximale $U_{F}$ de $-\Delta u+u^q=0$ dans un domaine quelconque de $\mathbb R^N$ impliquant la capacit\'e de Bessel $C_{2,q'}$ dans le cas sur-critique $q\geq q_{c}:=N/(N-2)$. Gr\^{a}ce \`a un crit\`ere de type Wiener, nous en d\'eduisons une condition n\'ecessaire et suffisante pour que cette solution maximale tende vers l'infini en un point du bord du domaine. Finalement nous prouvons un r\'esultat g\'en\'eral d'unicit\'e des grandes solutions.
{\it Pour citer cet article~: M. Marcus, L. V\'eron, C. R. Acad. Sci. Paris, Ser. I XXX (2006).}
\end{abstract}
\end{frontmatter}
\selectlanguage{francais}
\section*{Version fran\c{c}aise abr\'eg\'ee}
\setcounter {equation}{0}
Soit $F$ un sous-ensemble compact non-vide de  $\BBR^N$ de compl\'ementraire $F^c$ connexe et $q>1$. Il est bien connu qu'il existe une solution maximale $U_{F}$ de
\begin{equation}\label{Fq-eq}
 -\Gd u+u^q=0,
\end{equation}
dans $F^c=\BBR^N\setminus F$. En outre $U_{F}=0$ si et seulement si $C_{2,q'}\left(F\right)=0$, o\`u $q'=q/(q-1)$ et $C_{2,q'}$  d\'esigne la capacit\'e de Bessel en dimension $N$ \cite {BP}. Si $1<q<q_{c}:=N/(N-2)$, la capacit\'e de tout point est positive et la solution maximale est une grande solution \cite {Ve}, c'est \`a dire v\'erifie
\begin{equation}\label{expl}
\lim_{F^c\ni x\to y}U_{F}(x)=\infty,
\end{equation}
pour tout $y\in \prt F^c$, et la relation (\ref{expl}) est
uniforme en $y$. En outre $U_{F}$ est l'unique grande solution si
on suppose $\prt F^c\subset \prt \overline{F^c}$. Dans le cas
sur-critique $q\geq q_{c}$ la situation est beaucoup plus
compliqu\'ee dans la mesure o\`u les singularit\'es isol\'ees sont
\'eliminables et o\`u il existe une grande vari\'et\'e de
solutions. Si $q=2$, $N\geq 3$, Dhersin et Le Gall \cite{DL} ont
obtenu, par des m\'ethodes probabilistes, des estimations
pr\'ecises portant sur $U_{F}$ et utilisant la capacit\'e
$C_{2,2}$. De leurs estimations d\'ecoule une condition
n\'ecessaire et suffisante, exprim\'ee par un crit\`ere du type de
Wiener, pour que $U_{F}$ v\'erifie (\ref{expl}) en un point $y\in
\prt F^c$.

Labutin \cite{La} a r\'eussi \`a \'etendre partiellement les r\'esultats de \cite{DL} dans le cas $q\geq q_{c}$. Plus pr\'ecis\'ement il a prouv\'e que $U_{F}$ est une grande solution si et seulement si le crit\`ere de Wiener de
\cite{DL}, avec $C_{2,2}$ remplac\'e par $C_{2,q'}$, est v\'erifi\'e en {\it tout} point de $\prt F^c$, cependant il n'obtient pas l'estimation ponctuelle (\ref{expl}). Les estimations de Labutin sont optimales si $q>q_{c}$, mais pas si  $q=q_{c}$. Dans cette note nous \'etendons les r\'esultats de \cite {DL} par des m\'ethodes purement analytiques.

Si $F$ est un sous-ensemble ferm\'e non vide de $\BBR^N$, $x\in\BBR^N$ et $m\in\BBZ$ nous notons
$$T_{m}(x)=\left\{y\in\BBR^N:2^{-m-1}\leq\abs{x-y}\leq 2^{-m}\right\}$$ $$F_{m}(x)=F\cap T_{m}(x)\text { et }
 F^*_{m}(x)=F\cap \bar B_{2^{-m}}(x).$$
 On d\'efinit le {\it potentiel $C_{2,q'}$-capacitaire} $W_{F}$de $F$  par
 \begin{equation}\label{potcap}
W_{F}(x)=\sum_{-\infty}^\infty 2^\frac{2m}{q-1}C_{{2,q'}}(2^mF_{m}(x)).
\end{equation}
\noindent {\bf Theor\`eme 1.} {\it Il existe une constante $c=c(N,q)>0$
telle que}
 \begin{equation}\label{potcap1}
cW_{F}(x)\leq U_{F}(x)\leq \frac{1}{c}W_{F}(x)\forevery x\in F^c.
\end{equation}

Pour $q>q_{c}$ cette estimation est la m\^eme que celle de Labutin. Notre d\'emonstration est inspir\'ee de la sienne tout en faisant intervenir des arguments nouveaux qui simplifient notoirement sa d\'emarche. En utilisant la d\'efinition de la capacit\'e de Bessel on d\'emontre alors que la fonction $W_F$ est semi-continue sup\'erieurement dans $\overline {F^c}$. On en d\'eduit
\medskip

\noindent {\bf Theor\`eme 2.} {\it Pour tout point $y\in\prt F^c$,
 \begin{equation}\label{potcap2}
\lim_{F^c\ni x\to y}U_{F}(x)=\infty\Longleftrightarrow W_{F}(y)=\infty.
\end{equation}
Par suite $U_{F}$ est une grande solution si et seulement si $W_{F}(y)=\infty$ pour tout $y\in\prt F^c$.}\medskip

Il est facile de v\'erifier que si $W_{F}(y)=\infty$, alors $y$
est un point \'epais  de $F$, au sens de la topologie fine  $\frak
T_{q}$ associ\'ee \`a la capacit\'e $C_{{2,q'}}$. En utilisant la
propri\'et\'e de Kellog \cite{AH} que v\'erifie la capacit\'e
$C_{{2,q'}}$, on en d\'eduit que la solution maximale $U_{F}$ est
une {\it presque grande solution} dans le sens suivant: La
relation (\ref{expl}) a lieu sauf peut-\^etre sur un ensemble de
$\prt F^c$ de capacit\'e $C_{2,q'}$ nulle. \medskip

Il est classique que l'\'equation $-\Gd u+\abs u^{q-1}u=\gm$ admet une unique solution, not\'ee $u_{\gm}$, pour tout $\gm\in W^{-2,q'}(\BBR^N)$ \cite{BP}. On a alors le r\'esultat suivant\medskip

\noindent {\bf Theor\`eme 3.} {\it Pour tout sous-ensemble ferm\'e $F\subset\BBR^N$,
 \begin{equation}\label{potcap3}
U_{F}=\sup\{u_{\gm}:\gm\in W^{-2,q'}(\BBR^N),\gm(F^c)=0\}.
\end{equation}
Par suite $U_{F}$ est $\gs$-mod\'er\'ee, c'est \`a dire qu'il existe une suite croissante $\{\gm_{n}\}\subset W^{-2,q'}(\BBR^N)$ telle que
$\gm_{n}(F^c)=0$ et $u_{\gm_n}\uparrow u$.}\medskip

Cet \'enonc\'e est l'analogue dans le cas du probl\`eme elliptique
int\'erieur de  r\'esultats similaires concernant les probl\`emes
elliptique au bord \cite{MV3} et parabolique \cite {MV4}. Enfin,
 nous avons le r\'esultat d'unicit\'e suivant o\`u nous
d\'esignons par $\tilde E$ la fermeture de $E\subset\BBR^N$ pour
la topologie $\frak T_{q}$.\medskip

\noindent {\bf Theor\`eme 4.} {\it Pour tout ouvert non vide $D\subset\BBR^N$, posons $F=D^c$ et $F_{0}=\tilde D^c$ (c'est \`a dire que $F_{0}$ est l'int\'erieur de $F$ pour la topologie $\frak T_{q}$). Si
$C_{{2,q'}}(F\setminus \tilde F_{0})=0$, alors il existe au plus une grande solution de (\ref{Fq-eq}) dans $D$.}

\selectlanguage{english}
\section{Introduction}
\label{sec1}
\setcounter{equation}{0}
In this note we study  positive solutions of the equation
\begin{equation}\label{eq}
 -\Gd u+u^q=0,
\end{equation}
in $\BBR^N\setminus F$, $N\geq 3$, where $F$ is a non-empty compact set with $F^c$ connected and $q>1$. More precisely,
we shall study the behavior of the maximal solution of this problem, which we denote by $U_F$.
The existence of the maximal solution is guaranteed by the Keller-Osserman estimates (see \cite {MV2} for discussion about large solutions and the references therein).
\bcom
\begin{equation}\label{problem1}\BAL
-\Gd u+u^q&=0 &&\text{in}\q\Gw\sms F\\
u&=0 &&\text{on}\q\bdw
\EAL
\end{equation}
when $\Gw$ is either a smooth bounded domain or $\Gw=\BBR^N$.
This solution will be denoted by $U_F^\Gw$ or simply by $U_F$, if $\Gw=\BBR^N$. Note that $U_F^\Gw\leq U_F$.
\end{comment}
It is known \cite{BP} that, if $C_{2,q'}\left(F\right)=0$ then $U_F=0$.
If $u$ is a solution of \req{eq} in $D=\BBR^N\sms F$ and
$u$ blows up at every point of $\prt D$ we say that $u$ is a {\em large solution} in $D$. Obviously
a large solution exists in $D$ if and only if $U_F$ is a large solution.

 Our aim is:
(a) To provide a necessary and sufficient condition for the blow up of $U_F$  at an arbitrary  point $y\in F$ and
(b) to obtain a general uniqueness result for large solutions.

In the subcritical case, i.e. $1<q<q_c:=N/(N-2)$, these problems are well understood.
In this case $C_{2,q'}\left(F\right)>0$ for any non-empty set and
it is classical that positive solutions may have isolated point singularities of two
types: weak and strong.
This easily implies that the maximal solution $U_F$ is always a large solution in $D$. In addition it
is proved in   \cite{Ve} that the large solution is unique if it is assumed $\prt F^c\subset \prt \overline{F^c}^c$.


In the supercritical case, i.e. $q\geq q_{c}$, the situation is much more complicated.
 In this case  point singularities are removable and there exists a large variety of singular solutions.

 Sharp estimates for
 $U_F$ were obtained by Dhersin and Le Gall \cite{DL} in the case $q=2$, $N\geq 3$. These estimates
 were expressed in
terms of the Bessel capacity $C_{2,2}$
 and  were used to provide
 a Wiener type criterion for the pointwise blow up  of $U_F$, i.e., for  $y\in F$,
\begin{equation}\label{pt blowup}
 \lim_{F^c\ni x\to y} U_F(x)= \infty \iff \text{the Wiener type criterion is satisfied at y.} 
\end{equation}
These results were obtained by probabilistic tools; hence the restriction to $q=2$.

Labutin \cite{La} succeeded in partially extending the results of \cite{DL}
to $q\geq q_c$. Specifically, he proved that $U_F$ is a large
solution if and only if the Wiener  criterion of \cite{DL}, with
$C_{2,2}$ replaced by $C_{2,q'}$,  is satisfied \textit{at every
point of $F$}. The pointwise blow up was not established.
Labutin's result was obtained by analytic techniques. As in \cite{DL},
the proof is based on upper and lower estimates for $U_F$,  in
terms of the capacity $C_{2,q'}$. Labutin's estimates are sharp
for $q>q_c$ but not for $q= q_c$.

Conditions for uniqueness of  large solutions, for arbitrary $q>1$,  can be found in \cite{Ve}
and \cite{MV2}.

In the present paper we obtain a full extension of the results of \cite{DL}  to $q \geq q_c$,
$N\geq 3$.

Further we establish the following rather surprising fact: For any non-empty closed set $F\subsetneq \BBR^N$, the
maximal solution $U_F$ is
an 'almost large' solution in $D$ in the following sense: \req{pt blowup} holds at all points of $F$ with the
possible exception of a set
of $C_{2,q'}$-capacity zero. (Of course if $y$ is an interior point of $F$, \req{pt blowup}
holds in void.)

Finally we provide a capacitary sufficient condition for the uniqueness of large solutions.

\section{Statement of main results}
 Throughout the remainder of the note we assume that $q\geq q_{c}$. We start with some notation.
  For any set $A\sbs \BBR^N$ we denote by $\gr_A$ the distance function, $\gr_A(x)=\dist (x,A)$
   for every $x\in \BBR^N.$
If $F$ is a closed set and $x\in \BBR^N$ we denote
\begin{equation}\label{Fm} \BAL &T_m(x)=\set{y\in \BBR^N:2^{-(m+1)}\leq \abs{y-x}\leq 2^{-m}},\\
  &F_m(x)=F\cap T_m(x), \q F^*_m(x)=F\cap \bar B_{2^{-m}}(x).\EAL
\end{equation}
As usual $C_{\ga,p}$ denotes Bessel capacity in $\BBR^N$. Note that if $\ga=2$ and $p=q'=q/(q-1)$ then, for $q\geq N/(N-2)$,
$\ga p\leq N$. Put
\begin{equation}\label{WF0}
 W_F(x)= \sum_{-\infty}^{\infty} 2^{\frac{2m}{q-1}}C_{2,q'}\left(2^{m}F_m(x)\right).
\end{equation}
$W_F$ is called the {\em $C_{2,q'}$-capacitary potential} of $F$.

Observe that $2^{m}F^*_m(x)\sbs B_1(x)$ and that, for every $x\in F^c$, there exists a minimal integer
$M(x)$ such that $F_m(x)=\ems$ for $M(x)<m$. Therefore
\begin{equation}\label{WFa}
 W_F(x)= \sum_{-\infty}^{M(x)} 2^{\frac{2m}{q-1}}C_{2,q'}\left(2^{m}F_m(x)\right)<\infty \forevery x\in F^c
\end{equation}
It is known that there exists a constant $C$ depending only on $q,N$ such that
\begin{equation}\label{WF2}
 W_F(x)\leq W^*_F(x):=\sum_{m(x)}^{\infty}
2^{-\frac{2m}{q-1}}C_{2,q'}\left(2^{-m}F^{}*_m(x)\right)\leq C W_F(x)
\end{equation}
for every $x\in F^c$, see e.g. \cite{MV3}.

In the following results $F$ denotes a proper closed subset of $\BBR^N$.
The first theorem describes the capacitary estimates for the maximal solution.\medskip

\bth{capest} The maximal solution $U_F$ satisfies the inequalities
\begin{equation}\label{capest}
 \rec{c}W_F(x)\leq U_F(x)\leq cW_F(x)\forevery x\in F^c.
\end{equation}
\es
For $q>q_c$ these estimates are equivalent to those obtained by Labutin \cite{La}.
Our proof is inspired by the proof of \cite{La}, but employs some new arguments which lead to
a sharp estimate in the border case $q=q_c$ as well.
Using the previous theorem we establish:\medskip

\bth{large1} For every point $y\in F$,
\begin{equation}\label{Wiener}
 \lim_{F^c\ni x\to y} U_F(x)= \infty \iff W_F(y)=\infty.
\end{equation}
\Consy $U_F$ is a large solution in $F^c$ if and only if $W_F(y)=\infty$ for every $ y\in F.$
\es\medskip

\bth{q-ae} For any closed set $F\subsetneq\BBR^N$, the maximal solution $U_F$ is an almost large solution
in $D=F^c$ (see the definition of this term in the introduction).
\es
It is known \cite{BP} that if $\mu\in W^{-2,q}_+(\BBR^N)$ there exists  a unique solution of the equation
$
-\Gd u+u^q=\mu \txt{in} \BBR^N.
$
 This solution will be denoted by $u_\mu$.\medskip

\bth{U=V} For any closed set $F\subsetneq\BBR^N$,
\begin{equation}\label{mod-uF}
 U_F=\sup\set{u_\mu:\mu\in W^{-2,q}_+(\BBR^N),\; \mu(F^c)=0}.
\end{equation}
 Thus $U_F$ is \gsmod, i.e., there exists an increasing \seq $\set{\mu_n}\sbs W^{-2,q}_+(\BBR^N)$ such that
 $\mu_n(F^c)=0$
and $u_{\mu_n}\uar U_F$.
\es

For the next result we need the concept of
the $C_{2,q'}$-fine topology (in $\BBR^N$) that we shall denote by $\frak T_{q}$.  For its definition and basic properties see \cite[Ch. 6]{AH}.
The closure of a set $E$ in the
topology $\frak T_{q}$ will be denoted by $\tl E$.
The following uniqueness result holds.\medskip

\bth{unique}  Let $D\subset \BBR^N$ be
a non-empty, bounded open set. Put $F=D^c$ and $F_0=(\tl D)^c$ so that $F_0$ is the $\frak T_{q}$-interior of $F$.
If $C_{2,q'}\left(F\sms\tl F_0\right)=0$
then there exists at most one large solution in $D$.
\es

\section{Sketch of proofs.}
\noindent{\em On the proof of \rth{capest}.}\hskip 2mm The proof of this theorem is an adaptation of
the proof of the capacitary estimates for boundary value problems in \cite{MV3}. A central element of the proof in that paper
is the mapping $\BBP:W^{-2/q,q}_+(\bdw)\mapsto L^q(\Gw;\gr_{\bdw})$ given by
$\BBP(\mu)=\int_{\bdw}P(x,y)d\mu(y)$ where $P$ is the Poisson kernel in $\Gw$. In the  proof of the  present result the
same role is played by the Green operator acting on bounded measures in $\BBR^N$.\2
{\em On the proof of \rth{large1}.}\hskip 2mm
Denote
\begin{equation}\label{am(y)}
 a_m(x)=C_{2,q'}\left(2^{m}F_m(x)\right), \q a^*_m(x)=C_{2,q'}\left(2^{m}F^*_m(x)\right)
\end{equation}

First we show that $W_F(y)=\infty$ implies that $\lim_{D\ni x\to y}U_F(x)=\infty$.
Let $x\in F^c$ and let $\gl$ be an integer such that $2^{-\gl}\leq |x-y|\leq 2^{-\gl+1}$.
Obviously $\gl\leq M(x)$. For $m\leq \gl$:
 $$a^*_m(y)\leq C_{2,q'}\left(2(2^{m-1}{F^*}_m(x)\right)\leq ca^*_{m-1}(x).$$
Therefore
$$\sum_1^{\gl} 2^{\frac{2m}{q-1}}a^*_m(y)\leq
c\sum_1^{\gl} 2^{\frac{2m}{q-1}}a^*_{m-1}(x) \leq
c\sum_0^{\gl}2^{\frac{2m}{q-1}}a^*_{m}(x)\leq W^*_F(x).$$
As $x\to y$, $\gl \to\infty$ and the left hand side tends to $\infty$.
The reverse implication is a consequence of the following property of $W_F$.\medskip

\blemma{lcs} The function $y\mapsto W_F(y)$ is lower semi-continuous on $\overline{F^c}$. In addition, if   $W_F(y)<\infty$ then
$\liminf_{x\to y}W_F(x)<\infty$.
\es

\medskip

\noindent For proving this result, we use the fact that, for any $y\in\overline{F^c}$, and any $m\in\BBZ$,
$$C_{2,q'}\left(2^{m}F_m(y)\right)=
\inf\left\{\norm{\gz}^{q'}_{W^{2,q'}}:\gz\in C^\infty_0(\BBR^N):\gz\geq 0,\gz\geq 1\text { in a neighborhood of } 2^mF_m(y)\right\}.
$$
Thus, if $\gz\geq 1$ in a neighborhood of $2^mF_m(y)$, it implies that, for $\abs {x-y}$ small enough, $\gz\geq 1$ in a neighborhood of $2^mF_m(x)$. This implies
$$\lim_{\ge\to 0}\sup_{x\in\overline{F^c}\cap B_\ge(y)}
C_{2,q'}\left(2^{m}F_m(x)\right)=\limsup_{\overline{F^c}\ni x\to y}C_{2,q'}\left(2^{m}F_m(x)\right)\leq
C_{2,q'}\left(2^{m}F_m(y)\right).
$$
This implies the first assertion. The second assertion is proved by an argument involving the quasi-additivity
of capacity.\2
{\em On the proof of \rth{q-ae}.}\hskip 2mm It is not difficult to verify that, if $x$ is a thick point of
$F$ in the topology $\frak T_{q}$ (or $\frak T_{q}$-thick point),
then $W_F(x)=\infty$.
(For the definition of a thick point in a fine topology and the properties stated below see \cite[Ch. 6]{AH}.) The set of
$\frak T_{q}$-thick points of $F$ is denoted by $b_q(F)$ and it is known that, if $F$ is $\frak T_{q}$-closed
then $b_q(F)\sbs F$ and $C_{2,q'}\left(F\sms b_q(F\right)=0$ (this is called the Kellog property). Of course any set closed in the Euclidean topology is $\frak T_{q}$-closed.
Therefore, by \rth{capest}, $U_F$ blows up $C_{2,q'}$-a.e. on $\partial D$.\2
{\em On the proof of \rth{U=V}.}\hskip 2mm  Let us denote
the right hand side of \req{mod-uF} by $V_F$. Obviously $V_F\leq U_F$ and the proof of \rth{capest} actually shows that
$\rec{c}W_F(x)\leq V_F(x).$ Therefore $U_F\leq C V_F$ where $C$ is a constant depending only on $N,q$. By an argument introduced in \cite{MV1}
this implies that $U_F=V_F$.\2
{\em On the proof of \rth{unique}.}\hskip 2mm The proof is based on the following:\medskip

\blemma{UF=VF0} If $C_{2,q'}\left(F\sms\tl F_0\right)=0$ then
$$U_F=\sup\set{u_\mu:\mu\in W^{-2,q}_+(\BBR^N),\; \supp \mu\sbs F_0}.$$
\es The proof of the lemma involves
 subtle properties of the $C_{2,q'}$-fine topology.

\par The lemma implies that for every $x\in D$ there exists $\mu\in
W^{-2,q}_+(\BBR^N)$ such that $\supp \mu$ is a compact subset of
$F_0$ and $U_F(x)\leq Cu_\mu(x)$. Suppose that $u$ is a large
solution in $D$. Since $u_\mu$ is bounded in $\partial D$ it follows
that $u_\mu<u$. Thus $U_F\leq Cu$. By the argument of \cite{MV1}
mentioned before, this implies that $u=U_F$.

\vskip 2mm

\noindent{\bf Acknowledgment.} Both authors were
partially sponsored by an EC grant through the RTN Program
ñFront-Singularitiesî, HPRN-CT-2002-00274 and by the French-Israeli
cooperation program through grant No. 3-1352. The first author (MM)
also wishes to acknowledge the support of the Israeli Science
Foundation through grant No. 145-05.

\end{document}